# On the vertex-degree based invariants of digraphs*

Hanyuan Deng, Jiaxiang Yang, Zikai Tang, Jing Yang, Meiling You
Key Laboratory of Computing and Stochastic Mathematics (Ministry of Education),
College of Mathematics and Statistics, Hunan Normal University,
Changsha, Hunan 410081, P. R. China

**Abstract**

Let $D = (V, A)$ be a digraphs without isolated vertices. A vertex-degree based invariant $I(D)$ related to a real function $\varphi$ of $D$ is defined as a summation over all arcs, $I(D) = \frac{1}{2}\sum_{uv\in A}\varphi(d_u^+, d_v^-)$, where $d_u^+$ (resp. $d_u^-$) denotes the out-degree (resp. in-degree) of a vertex $u$. In this paper, we give the extremal values and extremal digraphs of $I(D)$ over all digraphs with $n$ non-isolated vertices. Applying these results, we obtain the extremal values of some vertex-degree based topological indices of digraphs, such as the Randić index, the Zagreb index, the sum-connectivity index, the $GA$ index, the $ABC$ index and the harmonic index, and the corresponding extremal digraphs.

**Keywords**: invariant; extremal digraph; Randić index; Zagreb index; sum-connectivity index; $GA$ index; $ABC$ index; harmonic index.

## 1 Introduction

A digraph $D = (V, A)$ is an ordered pair $(V, A)$ consisting of a non-empty finite set $V$ of vertices and a finite set $A$ of ordered pairs of distinct vertices called arcs (in particular, $D$ has no loops). If $a \in A$ is an arc from vertex $u$ to vertex $v$, then we indicate this by writing $a = uv$. The vertex $u$ is the tail of $a$, and the vertex $v$ its head. The out-degree (resp. in-degree) of a vertex $u$, denoted by $d_u^+$ (resp. $d_u^-$) is the number of arcs with tail $u$ (resp. with head $u$). A vertex $u$ for which $d_u^+ = d_u^- = 0$ is called an isolated vertex. We denote by $\mathcal{D}_n$ the set of all digraphs with $n$ non-isolated vertices.

Recently, J. Monsalve and J. Rada [7] extended the concept of vertex-degree based topological indices of graphs to digraphs. They obtained the extremal values of the Randić index of digraphs over $\mathcal{D}_n$, and found the extremal values of the Randić index over the set of all oriented trees with $n$ vertices. Also, they found the extremal values of the Randić index over the set of all orientations of the path, the cycle with $n$ vertices and the hypercube $H_d$ of dimension $d$, respectively.

All the digraphs considered in this paper are strict, i.e., no loops and no two arcs with the same ends have the same orientation.

A vertex-degree-based (VDB, for short) VDB invariant (or VDB topological index) $I(D)$ related to a real function $\varphi$ of a digraph $D$ with $n$ non-isolated vertices is defined as

$$I(D) = \frac{1}{2}\sum_{1\le i,j\le n-1} a_{ij}\varphi_{ij} \tag{1}$$

*This work is supported by the Hunan Provincial Natural Science Foundation of China (2020JJ4423), the Department of Education of Hunan Province (19A318) and the National Natural Science Foundation of China (11971164).

where $\varphi_{ij} = \varphi(i,j)$ and $a_{ij}$ is the number of arcs in $D$ of the form $uv$ such that $d_u^+ = i$ and $d_v^- = j$, i.e., $(i,j)$-arcs in $D$.

Recall that if $G$ is a graph, we can identify $G$ with the symmetric digraph $\overrightarrow{G}$ by replacing every edge of $G$ with a pair of symmetric arcs. Under this correspondence,

$$I(G) = \sum_{1\le i\le j\le n-1} m_{ij}\varphi_{ij} = I(\overrightarrow{G})$$

for any VDB topological index $\varphi$ with $\varphi_{ij} = \varphi_{ji}$ (symmetric) and $m_{ij}$ the number of edges in $G$ joining vertices of degree $i$ and $j$. In other words, The VDB topological index of digraphs is a generalization of the concept of VDB topological index of graphs.

In fact, a VDB topological index $I(D)$ of a digraph is an invariant based on the weights of all arcs depends on the out degrees of their tails and the in-degrees of their heads, i.e.,

$$I(G) = \sum_{uv\in A} \varphi(d_u^+, d_v^-)$$

where $\varphi(x,y)$ is a real function of $x$ and $y$ with $\varphi(x,y) \ge 0$ and $\varphi(x,y) = \varphi(y,x)$.

(i) If $\varphi(x,y) = (xy)^\alpha$, where $\alpha \ne 0$ is a real number, then $I(D)$ is the general Randić index of a digraph $D$. Furthermore, $I(D)$ is the Randić index, the second Zagreb index and the second modified Zagreb index for $\alpha = -\frac{1}{2}$, $\alpha = 1$ and $\alpha = -1$, respectively. For these indices of graphs, see [1, 6, 8, 9].

(ii) If $\varphi(x,y) = (x+y)^\alpha$, then $I(D)$ is the general sum-connectivity index of a digraph $D$. Further, $I(G)$ is the sum-connectivity index and the first Zagreb index for $\alpha = -\frac{1}{2}$ and $\alpha = 1$, respectively. See [5, 8, 11, 12] for graphs.

(iii) If $\varphi(x,y) = \frac{\sqrt{xy}}{\frac{1}{2}(x+y)}$, then $I(D)$ is the first geometric-arithmetic index $GA$ of a digraph $D$. See [10] for the first geometric-arithmetic index of a graph.

(iv) If $\varphi(x,y) = \sqrt{\frac{x+y-2}{xy}}$, then $I(G)$ is the atom-bond connectivity ($ABC$) index of a digraph $D$. See [3] for the atom-bond connectivity index of a graph.

(v) If $\varphi(x,y) = \frac{2}{x+y}$, then $I(D)$ is the harmonic index of a digraph $D$. See [4] for the harmonic index of a graph.

In this paper, we give the extremal values and extremal graphs of the VDB topological indices over all digraphs with $n$ non-isolated vertices by a unified linear-programming modeling, and provide a unified approach to determining some extremal values and characterizing extremal diraphs of Randić index, Zagreb index, sum-connectivity index, $GA$ index, $ABC$ index and harmonic index by using the linear programming methods.

## 2 General results on VDB invariants

Let $D$ be a digraph on $n \ge 2$ vertices without isolated vertices and $a_{ij}$ the number of arcs of $D$ from vertices of out-degree $i$ to vertices of in-degree $j$. If $\varphi$ is symmetric, i.e. $\varphi_{ij} = \varphi_{ji}$ for all $1 \le i < j \le n-1$, then we can simplify the expression in (1) in the following

$$I(D) = \frac{1}{2}\sum_{1\le i\le j\le n-1} p_{ij}\varphi_{ij} \tag{2}$$

where $p_{ij} = a_{ij} + a_{ji}$ for $i \neq j$ and $1 \leq i, j \leq n-1$, and $p_{ii} = a_{ii}$ for all $i = 1, 2, \cdots, n-1$.

Note that $p_{ij} = p_{ji}$ for all $1 \leq i, j \leq n-1$, and

$$\sum_{j=1}^{n-1} p_{ij} + p_{ii} = i n_i, \;\; 1 \leq i \leq n-1. \tag{3}$$

where $n_i$ is the number of vertices of $D$ with out-degree $i$ or in-degree $i$. Also,

$$n_1 + n_2 + \cdots + n_{n-1} = 2n - n_0. \tag{4}$$

The digraphs with $n$ non-isolated vertices which satisfy the following conditions are of great interest to us

(i)

$$\begin{cases} p_{ij} &= 0 \quad \text{for all } (i,j) \in \{(i,j) | 1 \leq i \leq j \leq n-1\} - \{(1, n-1)\}, \\ n_0 &= 0. \end{cases} \tag{5}$$

i.e., a digraph with only $(1, n-1)$- or $(n-1, 1)$-arcs and the out-degree or in-degree of each vertex greater than 0. The digraph obtained from the star on $n$ vertices by replacing each of its edges with a pair of symmetric arcs satisfies (5). The converse of this example does not hold since $D_1 = (V, A)$ is also a digraph satisfied (5), where $V = \{v_1, v_2, \cdots, v_n\}$ and $A = \{v_1v_2, v_2v_1, v_iv_1, v_2v_i | 3 \leq i \leq n\}$.

(ii)

$$\begin{cases} p_{ij} &= 0 \quad \text{for all } (i,j) \in \{(i,j) | 1 \leq i < j \leq n-1\}, \\ n_0 &= n. \end{cases} \tag{6}$$

i.e., the digraphs with only $(i,i)$-arcs $(1 \leq i \leq n-1)$ and the out-degree or in-degree of each vertex equal to 0. $\overrightarrow{K}_2$ satisfies (6), and $D_2 = (V, A)$ is also a digraph satisfied (6), where $V = \{v_1, v_2, v_3, v_4\}$ and $A = \{v_1v_3, v_1v_4, v_2v_3, v_2v_4\}$. All digraphs in which each component is $\overrightarrow{K}_2$ or $D_2$ satisfy (6).

(iii)

$$\begin{cases} p_{ij} &= 0 \quad \text{for all } (i,j) \in \{(i,j) | 1 \leq i < j \leq n-1\}, \\ n_0 &= 0. \end{cases} \tag{7}$$

i.e., the digraphs with only $(i,i)$-arcs $(1 \leq i \leq n-1)$ and the out-degree or in-degree of each vertex greater than 0. The dircted cycle $\overrightarrow{C}_n$ on $n$ vertices satisfies (7). All digraphs with $n$ non-isolated vertices in which each component is regular satisfy (7). The converse of this example does not hold since $D_3 = (V, A)$ is also a digraph satisfied (7), where $V = \{v_1, v_2, v_3, v_4\}$ and $A = \{v_1v_3, v_1v_4, v_2v_3, v_2v_4, v_3v_1, v_4v_2\}$.

We try to find $\min(I(G))$ and $\max(I(G))$ under the constraints (3) and (4). The following results give the solutions of this problem for some VDB topological indices $I(D)$, i.e., determine the extremal values and the correspond extremal digraphs of $I(D)$ over all digraphs on $n$ vertices without isolated vertices.

**Theorem 1.** *Let $D$ be a digraph on $n$ vertices without isolated vertices. $L_{ij} = \frac{n-1}{n}\left(\frac{1}{i} + \frac{1}{j}\right)\varphi_{1,n-1}$, and $S_1 = \{(i,j) | 1 \leq i \leq j \leq n-1\} - \{(1, n-1)\}$. Then*

*(i) If $\varphi_{ij} > L_{ij}$ for all $(i,j) \in S_1$, then $I(D) \geq \frac{n-1}{2}\varphi_{1,n-1}$ with equality if and only if $n_0 = n$ and $p_{ij} = 0$ for all $(i,j) \in S_1$, i.e., $D$ is the digraph $\overrightarrow{K}_{1,n-1}$ or $\overrightarrow{K}_{n-1,1}$, a star on $n$ vertices with its center of out-degree $n-1$ or 0.*

*(ii) If $\varphi_{ij} < L_{ij}$ for all $(i,j) \in S_1$, then $I(D) \leq (n-1)\varphi_{1,n-1}$ with equality if and only if $n_0 = 0$ and $p_{ij} = 0$ for all $(i,j) \in S_1$, i.e., $D$ satisfies the conditions (5).*

*Proof.* From (3), we have

$$n_i = \frac{1}{i}(\sum_{j=1}^{n-1} p_{ij} + p_{ii}), \quad i = 2, 3, \cdots, n-2, \tag{8}$$

$$n_1 - p_{1,n-1} = \sum_{j=1}^{n-2} p_{1j} + p_{11}, \tag{9}$$

$$(n-1)n_{n-1} - p_{1,n-1} = \sum_{j=2}^{n-1} p_{j,n-1} + p_{n-1,n-1}. \tag{10}$$

By (4) and (8),

$$n_1 + n_{n-1} = 2n - n_0 - \sum_{i=2}^{n-2} \frac{1}{i}(\sum_{j=1}^{n-1} p_{ij} + p_{ii}). \tag{11}$$

Multiplying (9) by $(n-1)$ and adding (10), we obtain

$$(n-1)(n_1 + n_{n-1}) - np_{1,n-1} = (n-1)\sum_{j=1}^{n-2} p_{1j} + (n-1)p_{11} + \sum_{j=2}^{n-1} p_{j,n-1} + p_{n-1,n-1},$$

and combining this equation with (11)

$$\begin{aligned} np_{1,n-1} = & \quad (n-1)(n_1 + n_{n-1}) - (n-1)(\sum_{j=1}^{n-2} p_{1j} + p_{11}) - (\sum_{j=2}^{n-1} p_{j,n-1} + p_{n-1,n-1}) \\ = & \quad (n-1)[2n - n_0 - \sum_{i=2}^{n-2} \tfrac{1}{i}(\sum_{j=1}^{n-1} p_{ij} + p_{ii})] - (n-1)(\sum_{j=1}^{n-2} p_{1j} + p_{11}) \\ & \quad -(\sum_{j=2}^{n-1} p_{j,n-1} + p_{n-1,n-1}). \end{aligned}$$

Hence,

$$\begin{aligned} p_{1,n-1} = & \quad 2(n-1) - \tfrac{n-1}{n} n_0 - \tfrac{n-1}{n} \sum_{i=2}^{n-2} \tfrac{1}{i}(\sum_{j=1}^{n-1} p_{ij} + p_{ii}) - \tfrac{n-1}{n}(\sum_{j=1}^{n-2} p_{1j} + p_{11}) \\ & \quad -\tfrac{1}{n} \sum_{j=2}^{n-1} p_{j,n-1} - \tfrac{1}{n} p_{n-1,n-1} \\ = & \quad 2(n-1) - \tfrac{n-1}{n} n_0 - \tfrac{n-1}{n}[\sum_{i=1}^{n-1} \tfrac{1}{i}(\sum_{j=1}^{n-1} p_{ij} + p_{ii}) - \tfrac{n}{n-1} p_{1,n-1}] \\ = & \quad 2(n-1) - \tfrac{n-1}{n} n_0 - \tfrac{n-1}{n}[\sum_{i=1}^{n-1} \tfrac{1}{i}(\sum_{j=1}^{n-1} p_{ij} + p_{ii})] + p_{1,n-1} \\ = & \quad 2(n-1) - \tfrac{n-1}{n} n_0 - \tfrac{n-1}{n}[\sum_{1 \le i \le j \le n-1} (\tfrac{1}{i} + \tfrac{1}{j}) p_{ij}] + p_{1,n-1} \\ = & \quad 2(n-1) - \tfrac{n-1}{n} n_0 - \tfrac{n-1}{n} \sum{}' (\tfrac{1}{i} + \tfrac{1}{j}) p_{ij} \end{aligned}$$

where $\sum'$ indicates summation over all $(i,j) \in S_1$. Substituting it into (2), we obtain

$$\begin{aligned} 2I(D) = & \quad \varphi_{1,n-1} p_{1,n-1} + \sum{}' \varphi_{ij} p_{ij} \\ = & \quad \varphi_{1,n-1}[2(n-1) - \tfrac{n-1}{n} n_0 - \tfrac{n-1}{n} \sum{}' (\tfrac{1}{i} + \tfrac{1}{j}) p_{ij}] + \sum{}' \varphi_{ij} p_{ij} \\ = & \quad [2(n-1) - \tfrac{n-1}{n} n_0] \varphi_{1,n-1} + \sum{}' [\varphi_{ij} - \tfrac{n-1}{n}(\tfrac{1}{i} + \tfrac{1}{j}) \varphi_{1,n-1}] p_{ij}. \end{aligned} \tag{12}$$

(i) If $\varphi_{ij} > L_{ij} = \frac{n-1}{n}(\frac{1}{i} + \frac{1}{j})\varphi_{1,n-1}$ for all $(i,j) \in S_1$, then (12) shows that $I(D) \geq \frac{1}{2}[2(n-1) - \frac{n-1}{n} n_0]\varphi_{1,n-1}$. Moreover, $I(D) \geq \frac{n-1}{2}\varphi_{1,n-1}$ since $n_0 \leq n$, with equality if and only if $n_0 = n$ and $p_{ij} = 0$ for all $(i,j) \in S_1$, i.e., $D$ is the digraph $\overrightarrow{K}_{1,n-1}$ or $\overrightarrow{K}_{n-1,1}$.

(ii) If $\varphi_{ij} < L_{ij}$ for all $(i,j) \in S_1$, then (12) shows that $I(D) \leq \frac{1}{2}[2(n-1) - \frac{n-1}{n} n_0]\varphi_{1,n-1}$. Moreover, $I(D) \leq (n-1)\varphi_{1,n-1}$ since $n_0 \geq 0$, with equality if and only if $n_0 = 0$ and $p_{ij} = 0$ for all $(i,j) \in S_1$, i.e., $D$ is a digraph satisfied (5). □

**Theorem 2.** *Let* $M_{ij} = \frac{n-1}{2}(\frac{1}{i}+\frac{1}{j})\varphi_{n-1,n-1}$, *and* $S_2 = \{(i,j)|1 \le i \le j \le n-1\} - \{(n-1,n-1)\}$. *Then*

*(i) If* $\varphi_{ij} > M_{ij}$ *for all* $(i,j) \in S_2$, *then* $I(D) \ge \frac{1}{4}n(n-1)\varphi_{n-1,n-1}$ *with equality if and only if* $n_0 = n$ *and* $p_{ij} = 0$ *for all* $(i,j) \in S_2$, *i.e.,* $D = \overrightarrow{K}_2$.

*(ii) If* $\varphi_{ij} < M_{ij}$ *for all* $(i,j) \in S_2$, *then* $I(D) \le \frac{1}{2}n(n-1)\varphi_{n-1,n-1}$ *with equality if and only if* $n_0 = 0$ *and* $p_{ij} = 0$ *for all* $(i,j) \in S_2$, *i.e.,* $D$ *is the digraph obtained from* $K_n$ *by replacing each edge with a pair of symmetric arcs.*

*Proof.* From (3) and (4), we obtain

$$\begin{aligned} n_{n-1} &= (2n-n_0) - \sum_{i=1}^{n-2}\tfrac{1}{i}(\sum_{j=1}^{n-1} p_{ij} + p_{ii}) \\ &= (2n-n_0) - (\sum_{j=1}^{n-1}\sum_{i=1}^{n-2}\tfrac{1}{i}p_{ij} + \sum_{i=1}^{n-2}\tfrac{1}{i}p_{ii}) \\ &= (2n-n_0) - \sum_{1\le i\le j\le n-1}(\tfrac{1}{i}+\tfrac{1}{j})p_{ij} + \tfrac{1}{n-1}\sum_{j=1}^{n-2}p_{j,n-1} + \tfrac{2}{n-1}p_{n-1,n-1}. \end{aligned}$$

By (3),

$$\sum_{j=1}^{n-2} p_{n-1,j} + 2p_{n-1,n-1} = (n-1)n_{n-1}$$

and

$$\begin{aligned} 2p_{n-1,n-1} &= (n-1)n_{n-1} - \sum_{j=1}^{n-2}p_{n-1,j} \\ &= (n-1)[(2n-n_0) - \sum_{1\le i\le j\le n-1}(\tfrac{1}{i}+\tfrac{1}{j})p_{ij} + \tfrac{1}{n-1}\sum_{j=1}^{n-2}p_{j,n-1} + \tfrac{2}{n-1}p_{n-1,n-1}] \\ &\quad - \sum_{j=1}^{n-2}p_{n-1,j} \\ &= (2n-n_0)(n-1) - (n-1)\sum\nolimits''(\tfrac{1}{i}+\tfrac{1}{j})p_{ij} \end{aligned}$$

where $\sum''$ indicates summation over all $(i,j) \in S_2$. By substituting it into (2), we obtain

$$\begin{aligned} 2I(D) &= \varphi_{n-1,n-1}p_{n-1,n-1} + \sum\nolimits''\varphi_{ij}p_{ij} \\ &= \varphi_{n-1,n-1}[\tfrac{1}{2}(2n-n_0)(n-1) - \tfrac{1}{2}(n-1)\sum\nolimits''(\tfrac{1}{i}+\tfrac{1}{j})p_{ij}] + \sum\nolimits''\varphi_{ij}p_{ij} \\ &= \tfrac{1}{2}(2n-n_0)(n-1)\varphi_{n-1,n-1} + \sum\nolimits''[\varphi_{ij} - \tfrac{n-1}{2}(\tfrac{1}{i}+\tfrac{1}{j})\varphi_{n-1,n-1}]p_{ij}. \end{aligned} \tag{13}$$

(i) If $\varphi_{ij} > M_{ij}$ for all $(i,j) \in S_2$, then (13) shows that $I(D) \ge \frac{1}{4}(2n-n_0)(n-1)\varphi_{n-1,n-1}$. Moreover, $I(D) \ge \frac{1}{4}n(n-1)\varphi_{n-1,n-1}$ since $n_0 \le n$, with equality if and only if $n_0 = n$ and $p_{ij} = 0$ for all $(i,j) \in S_2$, i.e., $D$ is a digraph with only $(n-1,n-1)$-arcs and the out-degree or in-degree of each vertex equal to 0. So, $D = \overrightarrow{K}_2$.

(ii) If $\varphi_{ij} < M_{ij}$ for all $(i,j) \in S_2$, then (13) shows that $I(D) \le \frac{1}{4}(2n-n_0)(n-1)\varphi_{n-1,n-1}$. Moreover, $I(D) \le \frac{1}{2}n(n-1)\varphi_{n-1,n-1}$ since $n_0 \ge 0$, with equality if and only if $n_0 = 0$ and $p_{ij} = 0$ for all $(i,j) \in S_2$, i.e., $D$ is the digraph obtained from the complete graph $K_n$ by replacing each edge with a pair of symmetric arcs. □

**Theorem 3.** *Let* $S_3 = \{(i,j)|1 \le i \ne j \le n-1\}$. *Then*

*(i) If $\varphi_{ij} > M_{ij}$ for all $(i,j) \in S_3$, and $i\varphi_{ii} = (n-1)\varphi_{n-1,n-1}$ for $1 \le i \le n-2$, then $I(D) \ge \frac{1}{4}n(n-1)\varphi_{n-1,n-1}$ with equality if and only if $n_0 = n$ and $p_{ij} = 0$ for all $(i,j) \in S_3$, i.e., $D$ is a digraph satisfied (6).*

*(ii) If $\varphi_{ij} < M_{ij}$ for all $(i,j) \in S_3$, and $i\varphi_{ii} = (n-1)\varphi_{n-1,n-1}$ for $1 \le i \le n-2$, then $I(D) \le \frac{1}{2}n(n-1)\varphi_{n-1,n-1}$ with equality if and only if $n_0 = 0$ and $p_{ij} = 0$ for all $(i,j) \in S_3$, i.e., $D$ is a digraph satisfied (7).*

*Proof.* From (13), we have

$$\begin{aligned} 2I(G) = \ & \tfrac{1}{2}(2n-n_0)(n-1)\varphi_{n-1,n-1} + \sum_{1\le i<j\le n-1}[\varphi_{ij} - \tfrac{n-1}{2}(\tfrac{1}{i}+\tfrac{1}{j})\varphi_{n-1,n-1}]p_{ij} \\ & + \sum_{i=1}^{n-2}[\varphi_{ii} - \tfrac{n-1}{i}\varphi_{n-1,n-1}]p_{ii}. \end{aligned} \tag{14}$$

(i) If $\varphi_{ij} > M_{ij}$ for all $(i,j) \in S_3$, and $i\varphi_{ii} = (n-1)\varphi_{n-1,n-1}$ for $1 \le i \le n-2$, then (14) shows that $I(D) \ge \frac{1}{4}(2n-n_0)(n-1)\varphi_{n-1,n-1}$. Moreover, $I(D) \ge \frac{1}{4}n(n-1)\varphi_{n-1,n-1}$ since $n_0 \le n$, with equality if and only if $n_0 = n$ and $p_{ij} = 0$ for all $(i,j) \in S_3$, i.e., $D$ is the digraph satisfied (6).

(ii) If $\varphi_{ij} < M_{ij}$ for all $(i,j) \in S_3$, and $i\varphi_{ii} = (n-1)\varphi_{n-1,n-1}$ for $1 \le i \le n-2$, then (14) shows that $I(D) \le \frac{1}{4}(2n-n_0)(n-1)\varphi_{n-1,n-1}$. Moreover, $I(D) \le \frac{1}{2}n(n-1)\varphi_{n-1,n-1}$ since $n_0 \ge 0$, with equality if and only if $n_0 = 0$ and $p_{ij} = 0$ for all $(i,j) \in S_3$, i.e., $D$ is a digraph satisfied (7). □

Theorems 1-3 show that the results on digraphs are different from the results on graphs in [2].

# 3 Applications

In this section, we give some results on Randić index, Zagreb index, sum-connectivity index, $GA$ index and $ABC$ index of digraphs by using Theorems 1-3.

## 3.1 The general Randić index of digraphs

Let $\varphi_{ij} = (ij)^\alpha$. Then $I(D) = R_\alpha(D) = \frac{1}{2}\sum_{1\le i\le j\le n-1} p_{ij}(ij)^\alpha$ is the general Randić index of a digraph $D$ with $n$ non-isolated vertices. In particular, $R_\alpha(D)$ is the Randić index, the second Zagreb index and the modified Zagreb index of a digraph for $\alpha = -\frac{1}{2}$, $\alpha = 1$ and $\alpha = -1$, respectively.

(i) Let $-\frac{1}{2} \le \alpha < +\infty$. Then $2\alpha + 1 \ge 0$.

Note that $ij \le (\frac{i+j}{2})^2$ and $i,j \le n-1$, we have

$$\begin{aligned} (ij)^{\alpha+1} & \le (\tfrac{i+j}{2})^{2\alpha+2} = \tfrac{1}{2^{2\alpha+2}}(i+j)^{2\alpha+1}(i+j) \\ & \le \tfrac{1}{2^{2\alpha+2}}[2(n-1)]^{2\alpha+1}(i+j) = \tfrac{1}{2}(n-1)^{2\alpha+1}(i+j), \end{aligned}$$

and $\varphi_{ij} \le \frac{n-1}{2}(\frac{1}{i}+\frac{1}{j})\varphi_{n-1,n-1}$ with equality if and only if (a) $i = j = n-1$ for $-\frac{1}{2} < \alpha < +\infty$, or (b) $i = j$ for $\alpha = -\frac{1}{2}$.

By Theorems 2(ii) and 3(ii), we have

$$R_\alpha(D) \le \frac{1}{2}n(n-1)\varphi_{n-1,n-1} = \frac{1}{2}n(n-1)^{2\alpha+1}$$

with equality if and only if (a) $D$ is the digraph obtained from $K_n$ by replacing each edge with a pair of symmetric arcs for $-\frac{1}{2} < \alpha < +\infty$, or (b) $D$ is a digraph satisfied (7).

So, (a) the digraph with the maximal general Randić index (including the second Zagreb index) for $-\frac{1}{2} < \alpha < +\infty$ is the digraph obtained from $K_n$ by replacing each edge with a pair of symmetric arcs; (b) the digraphs with the maximal Randić index are those satisfied (7), see Theorem 3.7 in [7].

**Corollary 4.** *If $D \in \mathcal{D}_n$, then (a) $R_\alpha(D) \le \frac{1}{2}n(n-1)^{2\alpha+1}$ for $-\frac{1}{2} < \alpha < +\infty$ with equality if and only if $D$ is the digraph obtained from $K_n$ by replacing each edge with a pair of symmetric arcs; (b) (Theorem 3.7 in [7]) $R_{-\frac{1}{2}}(D) \le \frac{n}{2}$ with equality if and only if $D$ satisfies (7).*

(ii) Let $-\infty < \alpha \le -1$.

Because $ij \le (\frac{i+j}{2})^2$ and $\alpha \le -1$, we have

$$\begin{aligned}(ij)^{\alpha+1} &\ge (\tfrac{i+j}{2})^{2\alpha+2} = \tfrac{1}{2^{2\alpha+2}}(i+j)^{2\alpha+1}(i+j)\\ &\ge \tfrac{1}{2^{2\alpha+2}}[2(n-1)]^{2\alpha+1}(i+j) = \tfrac{1}{2}(n-1)^{2\alpha+1}(i+j),\end{aligned}$$

and $\varphi_{ij} \ge \frac{n-1}{2}(\frac{1}{i}+\frac{1}{j})\varphi_{n-1,n-1}$ with equality if and only if $i = j = n-1$.

By Theorem 2(i), we have

$$R_\alpha(D) \ge \frac{1}{4}n(n-1)\varphi_{n-1,n-1} = \frac{1}{4}n(n-1)^{2\alpha+1}$$

with equality if and only if $D = \overrightarrow{K}_2$.

So, the digraph with the minimal general Randić index (including the modified Zagreb index) for $-\infty < \alpha \le -1$ is $\overrightarrow{K}_2$.

**Corollary 5.** *If $D \in \mathcal{D}_n$, then $R_\alpha(D) \le \frac{1}{4}n(n-1)^{2\alpha+1}$ for $-\infty < \alpha \le -1$ with equality if and only if $D = \overrightarrow{K}_2$.*

(iii) Let $-\frac{1}{2} \le \alpha < 0$.

In the following, we show that $\varphi_{ij} > \frac{n-1}{n}(\frac{1}{i}+\frac{1}{j})\varphi_{1,n-1}$ for all $(i,j) \in \{(i,j)|1 \le i \le j \le n-1\} - \{(1,n-1)\}$.

Let $g(x,y) = \frac{(xy)^{\alpha+1}}{x+y}$, where $1 \le x \le y \le n-1$. Note that $\alpha x + y + \alpha y \ge (2\alpha+1)x \ge 0$, $\frac{\partial g}{\partial x} = \frac{y(xy)^\alpha(\alpha x+y+\alpha y)}{(x+y)^2} = 0$ and $\frac{\partial g}{\partial y} = \frac{x(xy)^\alpha(\alpha x+x+\alpha y)}{(x+y)^2} = 0$ if and only if $\alpha = -\frac{1}{2}$ and $x = y$. So, the minimal point of $g(x,y)$ in the region $\{(x,y)|1 \le x \le y \le n-1\}$ is on the boundary of this region, and the minimal value of $g(x,y)$ in the region $\{(x,y)|1 \le x \le y \le n-1\}$ is $\min\{g(1,1), g(1,n-1)\} = \min\{\frac{1}{2}, \frac{(n-1)^{\alpha+1}}{n}\}$.

If $\alpha \in (-\frac{1}{2}, 0)$, then $\frac{(n-1)^{\alpha+1}}{n} < \frac{1}{2}$ for sufficiently large $n$; and if $\alpha = -\frac{1}{2}$, then $\frac{(n-1)^{\alpha+1}}{n} < \frac{1}{2}$ for $n \ge 3$. Hence, $g(i,j) \ge g(1,n-1)$, and

$$(ij)^\alpha \ge \frac{(n-1)^{\alpha+1}}{n}(\frac{1}{i}+\frac{1}{j}), \text{ i.e. } \varphi_{ij} \ge \frac{n-1}{n}(\frac{1}{i}+\frac{1}{j})\varphi_{1,n-1}$$

with equality if and only if $(i,j) = (1,n-1)$.

By Theorem 1(i), we have

$$R_\alpha(D) \ge \frac{n-1}{2}\varphi_{1,n-1} = \frac{1}{2}(n-1)^{\alpha+1}$$

with equality if and only if $D$ is the digraph $\overrightarrow{K}_{1,n-1}$ or $\overrightarrow{K}_{n-1,1}$ for sufficiently large $n$.

So, the digraph with the minimal Randić index is $\overrightarrow{K}_{1,n-1}$ or $\overrightarrow{K}_{n-1,1}$ over $\mathcal{D}_n$ for $n \geq 3$; and the digraph with the minimal general Randić index for $\alpha \in (-\frac{1}{2}, 0)$ is also $\overrightarrow{K}_{1,n-1}$ or $\overrightarrow{K}_{n-1,1}$ over $\mathcal{D}_n$ for sufficiently large $n$.

**Corollary 6.** *(a) (Theorem 3.11 in [7]) If $D \in \mathcal{D}_n$, $n \geq 3$, then $R_{-\frac{1}{2}}(D) \geq \frac{1}{2}\sqrt{n-1}$ with equality if and only if $D = \overrightarrow{K}_{1,n-1}$ or $D = \overrightarrow{K}_{n-1,1}$;*

*(b) Let $-\frac{1}{2} \leq \alpha < 0$. If $D \in \mathcal{D}_n$, then $R_\alpha(D) \geq \frac{1}{2}(n-1)^{\alpha+1}$ for sufficiently large $n$, with equality if and only if $D = \overrightarrow{K}_{1,n-1}$ or $D = \overrightarrow{K}_{n-1,1}$.*

## 3.2 The general sum-connectivity index of digraphs

Let $\varphi_{ij} = (i+j)^\alpha$. Then $I(D) = \chi_\alpha(D) = \frac{1}{2} \sum\limits_{1 \leq i \leq j \leq n-1} p_{ij}(i+j)^\alpha$ is the general sum-connectivity index of a digraph $D$, and $\chi_\alpha(D)$ is the sum-connectivity index and the first Zagreb index of $D$ for $\alpha = -\frac{1}{2}$ and $\alpha = 1$, respectively.

(i) Let $-1 \leq \alpha < +\infty$.

Because $1 \leq i \leq j \leq n-1$ and $\alpha + 1 \geq 0$,

$$\begin{array}{rl} ij & \leq (\frac{i+j}{2})^2 = (\frac{i+j}{2})^{1-\alpha}(\frac{i+j}{2})^{1+\alpha} \\ & \leq (\frac{i+j}{2})^{1-\alpha}(n-1)^{1+\alpha}, \end{array}$$

and $\varphi_{ij} = (i+j)^\alpha \leq \frac{n-1}{2}(\frac{1}{i} + \frac{1}{j})[2(n-1)]^\alpha = \frac{n-1}{2}(\frac{1}{i} + \frac{1}{j})\varphi_{n-1,n-1}$ with equality if and only if (a) $i = j = n-1$ for $-1 < \alpha < +\infty$, or (b) $i = j$ for $\alpha = -1$.

By Theorems 2(ii) and 3(ii), we have

$$\chi_\alpha(D) \leq \frac{1}{2}n(n-1)\varphi_{n-1,n-1} = 2^{\alpha-1}n(n-1)^{\alpha+1}$$

with equality if and only if (a) $D$ is the digraph obtained from the complete graph $K_n$ by replacing each edge with a pair of symmetric arcs, or (b) $D$ satisfies (7).

Especially, this shows that the graph with the maximal sum-connectivity index, or the maximal first Zagreb index is $K_n$ among all graphs of order $n$.

**Corollary 7.** *If $D \in \mathcal{D}_n$, then (a) $\chi_\alpha(D) \leq 2^{\alpha-1}n(n-1)^{\alpha+1}$ for $-\frac{1}{2} < \alpha < +\infty$ with equality if and only if $D$ is the digraph obtained from $K_n$ by replacing each edge with a pair of symmetric arcs; (b) $\chi_{-1}(D) \leq \frac{n}{4}$ with equality if and only if $D$ satisfies (7).*

(ii) Let $-1 \leq \alpha < 0$.

We consider the function $g(x, y) = (xy)(x+y)^{\alpha-1}$, where $1 \leq x \leq y \leq n-1$. It is easy to know that the minimal value of $g(x, y) = (xy)(x+y)^{\alpha-1}$ in the region $\{(x, y)|1 \leq x \leq y \leq n-1\}$ is $\min\{g(1, 1), g(1, n-1)\} = \min\{2^{\alpha-1}, (n-1)n^{\alpha-1}\}$.

If $\alpha \in (-\frac{1}{2}, 0)$, then $(n-1)n^{\alpha-1} < 2^{\alpha-1}$ for sufficiently large $n$; and if $\alpha \in [-1, -\frac{1}{2}]$, then $(n-1)n^{\alpha-1} < 2^{\alpha-1}$ for $n \geq 6$. Hence, $g(i, j) \geq g(1, n-1)$, and

$$(i+j)^\alpha \geq (n-1)n^{\alpha-1}(\frac{1}{i} + \frac{1}{j}), \text{ i.e., } \quad \varphi_{ij} \geq \frac{n-1}{n}(\frac{1}{i} + \frac{1}{j})\varphi_{1,n-1}$$

with equality if and only if $(i,j)=(1,n-1)$. By Theorem 1(i), we have

$$\chi_\alpha(D) \geq \frac{1}{2}(n-1)\varphi_{1,n-1} = \frac{1}{2}(n-1)n^\alpha$$

with equality if and only if $D$ is $\overrightarrow{K}_{1,n-1}$ or $\overrightarrow{K}_{n-1,1}$ for $\alpha \in [-1,0)$ and sufficiently large $n$, or for $\alpha \in [-1,-\frac{1}{2}]$ and $n \geq 6$.

So, that the graph with the minimal general sum-connectivity index for $\alpha \in [-1,-\frac{1}{2}]$ is $\overrightarrow{K}_{1,n-1}$ or $\overrightarrow{K}_{n-1,1}$ over $\mathcal{D}_n$; and the digraph with the minimal general sum-connectivity index for $\alpha \in [-1,0)$ is also $\overrightarrow{K}_{1,n-1}$ or $\overrightarrow{K}_{n-1,1}$ over $\mathcal{D}_n$ when $n$ is sufficiently large.

**Corollary 8.** *Let $D \in \mathcal{D}_n$. If $\alpha \in [-1,-\frac{1}{2}]$ and $n \geq 6$, or $\alpha \in [-1,0)$ and $n$ is sufficiently large, then $\chi_\alpha(D) \geq \frac{1}{2}(n-1)n^\alpha$ with equality if and only if $D$ is $\overrightarrow{K}_{1,n-1}$ or $\overrightarrow{K}_{n-1,1}$.*

## 3.3 The geometric-arithmetic index of digraphs

Let $\varphi_{ij} = \frac{\sqrt{ij}}{\frac{1}{2}(i+j)}$. Then $I(D) = GA(D) = \frac{1}{2}\sum\limits_{1\leq i\leq j\leq n-1} p_{ij}\frac{\sqrt{ij}}{\frac{1}{2}(i+j)}$ is the first geometric-arithmetic index $GA$ of a digraph $D$.

(i) Note that $\varphi_{n-1,n-1}=1$ and $(ij)^{\frac{3}{2}} \leq (\frac{i+j}{2})^3 = \frac{i+j}{8}(i+j)^2 \leq \frac{n-1}{4}(i+j)^2$, i.e. $\frac{\sqrt{ij}}{\frac{1}{2}(i+j)} \leq \frac{n-1}{2}(\frac{1}{i}+\frac{1}{j})$, we have $\varphi_{ij} \leq \frac{n-1}{2}(\frac{1}{i}+\frac{1}{j})\varphi_{n-1,n-1}$ with equality if and only if $i=j=n-1$. By Theorem 2(ii),

$$GA(D) \leq \frac{1}{2}n(n-1)\varphi_{n-1,n-1} = \frac{1}{2}n(n-1)$$

with equality if and only if $D$ is the digraph obtained from $K_n$ by replacing each edge with a pair of symmetric arcs.

(ii) It is easy to know that the minimal value of $g(x,y) = \frac{(xy)^{\frac{3}{2}}}{(x+y)^2}$ in the region $\{(x,y)|1 \leq x \leq y \leq n-1\}$ is $g(1,n-1) = \frac{(n-1)^{\frac{3}{2}}}{n^2}$, $g(i,j) \geq g(1,n-1)$, i.e. $\frac{(ij)^{\frac{3}{2}}}{(i+j)^2} \geq \frac{(n-1)^{\frac{3}{2}}}{n^2}$. Hence,

$$\frac{\sqrt{ij}}{\frac{1}{2}(i+j)} \geq \frac{n-1}{n}(\frac{1}{i}+\frac{1}{j})\frac{\sqrt{n-1}}{\frac{1}{2}n} \quad \text{i.e.,} \quad \varphi_{ij} \geq \frac{n-1}{n}(\frac{1}{i}+\frac{1}{j})\varphi_{1,n-1}$$

with equality if and only if $(i,j)=(1,n-1)$.

By Theorem 1(i), we have

$$GA(D) \geq \frac{1}{2}(n-1)\varphi_{1,n-1} = \frac{(n-1)^{\frac{3}{2}}}{n}$$

with equality if and only if $D$ is $\overrightarrow{K}_{1,n-1}$ or $\overrightarrow{K}_{n-1,1}$.

So, we obtain the digraphs with the maximal and the minimal geometric-arithmetic index $GA$ over $\mathcal{D}_n$.

**Corollary 9.** *If $D \in \mathcal{D}_n$, then $GA(D) \leq \frac{1}{2}n(n-1)$ with equality if and only if $D$ is the digraph obtained from $K_n$ by replacing each edge with a pair of symmetric arcs; $GA(D) \geq \frac{(n-1)^{\frac{3}{2}}}{n}$ with equality if and only if $D$ is $\overrightarrow{K}_{1,n-1}$ or $\overrightarrow{K}_{n-1,1}$.*

### 3.4 The atom-bond connectivity index of digraphs

Let $\varphi_{ij} = \sqrt{\frac{i+j-2}{ij}}$, then $I(D) = ABC(D) = \frac{1}{2} \sum\limits_{1\le i\le j\le n-1} p_{ij}\sqrt{\frac{i+j-2}{ij}}$ is the ABC index of a digraph $D$. Since $1 \le i \le j \le n-1$,

$$\frac{i+j-2}{ij} \le \frac{2(n-2)}{ij} \le \frac{2(n-2)}{(ij)^2}(\frac{i+j}{2})^2 \le \frac{n-2}{2}(\frac{i+j}{ij})^2,$$

and $\sqrt{\frac{i+j-2}{ij}} \le \sqrt{\frac{n-2}{2}}(\frac{i+j}{ij})$, i.e., $\varphi_{ij} \le \frac{n-1}{2}(\frac{1}{i}+\frac{1}{j})\varphi_{n-1,n-1}$ with equality if and only if $i = j = n-1$. By Theorem 2(ii), we have

$$ABC(D) \le \frac{1}{2}n(n-1)\varphi_{n-1,n-1} = \frac{1}{2}n\sqrt{2n-4}$$

with equality if and only if $D$ is the digraph obtained from $K_n$ by replacing each edge with a pair of symmetric arcs.

This shows that the digraphs with the maximal ABC index over $\mathcal{D}_n$ is the digraph obtained from $K_n$ by replacing each edge with a pair of symmetric arcs.

**Corollary 10.** *If $D \in \mathcal{D}_n$, then $ABC(D) \le \frac{1}{2}n\sqrt{2n-4}$ with equality if and only if $D$ is the digraph obtained from $K_n$ by replacing each edge with a pair of symmetric arcs.*

### 3.5 The harmonic index of digraphs

Let $\varphi_{ij} = \frac{2}{i+j}$. Then $I(D) = h(D) = \frac{1}{2} \sum\limits_{1\le i\le j\le n-1} p_{ij}\frac{2}{i+j}$ is the harmonic index of a digraph $D$.

(i) Note that $\varphi_{ij} = \frac{2}{i+j} \le \frac{i+j}{2ij} = \frac{n-1}{2}(\frac{1}{i}+\frac{1}{j})\varphi_{n-1,n-1}$ with equality if and only if $i = j$, and $i\varphi_{ii} = 1 = (n-1)\varphi_{n-1,n-1}$, from Theorem 3(ii), we have

$$h(D) \le \frac{1}{2}n(n-1)\varphi_{n-1,n-1} = \frac{n}{2}$$

with equality if and only if $D$ is a digraph satisfied (7).

(ii) Also, the minimal value of $g(x, y) = \frac{(xy)}{(x+y)^2}$ in the region $\{(x, y)|1 \le x \le y \le n-1\}$ is $g(1, n-1) = \frac{n-1}{n^2}$, we have $g(i, j) = \frac{ij}{(i+j)^2} \ge \frac{n-1}{n^2} = g(1, n-1)$. And $\varphi_{ij} = \frac{2}{i+j} \ge \frac{n-1}{n}(\frac{1}{i}+\frac{1}{j})\frac{2}{n} = \frac{n-1}{n}(\frac{1}{i}+\frac{1}{j})\varphi_{1,n-1}$ with equality if and only if $(i, j) = (1, n-1)$.

By Theorem 1(i), we have

$$h(D) \ge \frac{1}{2}(n-1)\varphi_{1,n-1} = \frac{n-1}{n}$$

with equality if and only if $D$ is $\overrightarrow{K}_{1,n-1}$ or $\overrightarrow{K}_{n-1,1}$.

So, we obtain the digraphs with the minimal and maximal harmonic index over $\mathcal{D}_n$.

**Corollary 11.** *If $D \in \mathcal{D}_n$, then $h(D) \le \frac{n}{2}$ with equality if and only if $D$ is a digraph satisfied (7); $h(D) \ge \frac{n-1}{n}$ with equality if and only if $D$ is $\overrightarrow{K}_{1,n-1}$ or $\overrightarrow{K}_{n-1,1}$.*